\documentclass[letterpaper, 10pt, conference]{IEEEconf}
\IEEEoverridecommandlockouts 

\usepackage{amsmath,amssymb}
\usepackage{epsfig,curves,multirow,subfigure,url}

\newcommand{\norm}[1]{\ensuremath{\left\| #1 \right\|}}
\newcommand{\bracket}[1]{\ensuremath{\left[ #1 \right]}}
\newcommand{\braces}[1]{\ensuremath{\left\{ #1 \right\}}}
\newcommand{\parenth}[1]{\ensuremath{\left( #1 \right)}}
\newcommand{\refeqn}[1]{(\ref{eqn:#1})}

\newcommand{\deriv}[2]{\ensuremath{\frac{\partial #1}{\partial #2}}}
\newcommand{\T}{\ensuremath{\mathrm{T}}}
\newcommand{\SO}{\ensuremath{\mathrm{SO(3)}}}
\newcommand{\so}{\ensuremath{\mathfrak{so}(3)}}
\newcommand{\SE}{\ensuremath{\mathrm{SE(3)}}}
\newcommand{\se}{\ensuremath{\mathfrak{se}(3)}}

\newcommand{\aSE}[2]{\ensuremath{\begin{bmatrix}#1&#2\\0&1\end{bmatrix}}}
\newcommand{\ase}[2]{\ensuremath{\begin{bmatrix}#1&#2\\0&0\end{bmatrix}}}

\renewcommand{\Re}{\ensuremath{\mathbb{R}}}
\renewcommand{\S}{\ensuremath{\mathbb{S}}}

\renewcommand{\theenumi}{\roman{enumi}}
\renewcommand{\labelenumi}{(\theenumi)}

\title{\LARGE \bf
Optimal Control of a Rigid Body\\ using Geometrically Exact
Computations on \SE}

\author{ \parbox{3 in}{\centering Taeyoung Lee\authorrefmark{1}\authorrefmark{2}, N. Harris McClamroch\authorrefmark{2}\\
         Department of Aerospace Engineering\\
         University of Michigan, Ann Arbor, MI 48109\\
         {\tt\small \{tylee, nhm\}@umich.edu}}
         \hspace*{ 0.5 in}
         \parbox{3 in}{\centering Melvin Leok\authorrefmark{1}\\
         Department of Mathematics\\
        Purdue University, West Lafayette, IN 47907\\
         {\tt\small mleok@math.purdue.edu}}
        \thanks{\textsuperscript{\footnotesize\ensuremath{*}}This research has been supported in part by NSF under grant DMS-0504747, and by a grant from the Rackham Graduate School, University of Michigan.}
        \thanks{\textsuperscript{\footnotesize\ensuremath{\dagger}}This research has been supported in part by NSF under grant ECS-0244977.}
}

\begin{document}
\maketitle \thispagestyle{empty} \pagestyle{empty}

\begin{abstract}
Optimal control problems are formulated and efficient computational
procedures are proposed for combined orbital and rotational
maneuvers of a rigid body in three dimensions.   The rigid body is
assumed to act under the influence of forces and moments that arise
from a potential and from control forces and moments.  The key
features of this paper are its use of computational procedures that
are guaranteed to preserve the geometry of the optimal solutions.
The theoretical basis for the computational procedures is
summarized, and examples of optimal spacecraft maneuvers are
presented.
\end{abstract}

\section{Introduction}
Discrete optimal control problems for translational and rotational
dynamics of a rigid body under a potential are studied. Optimal
control of a rigid body arises in numerous engineering and
scientific fields. These problems provide both a theoretical
challenge and a numerical challenge in the sense that the
configuration space has a Lie group structure denoted by $\SE$ that
defines a fundamental constraint.

Optimal control problems on a Lie group have been studied
in~\cite{Spin.MCSS98,Sas.ICIAM95}. These studies are based on the
driftless kinematics of a Lie group. The dynamics are ignored, and
it is assumed that elements in the corresponding Lie algebra are
controlled directly.

General-purpose numerical integration methods, including the popular
Runge--Kutta schemes, typically preserve neither the group structure
of the configuration space nor geometric invariants of the dynamics.
Geometric structure-preserving integrators, referred to as Lie group
variational integrators~\cite{pro:cca05}, preserve the group
structure without the use of local charts, reprojection, or
constraints, and they have the desirable property that they are
symplectic and momentum preserving, and they exhibit good energy
behavior for an exponentially long time period.

This paper presents geometrically exact and numerically efficient
computational approaches to solve optimal control problems of a
rigid body on a Lie group, $\SE$. The dynamics and the kinematics
are discretized by
 a Lie group variational integrator, and discrete
 optimality conditions are constructed. Efficient numerical
 algorithms to solve the necessary condition are developed.
This method provide a substantial advantage over current methods for
optimal control on a Lie group in the sense that the dynamics of a
rigid body as well as the kinematics equation are explicitly
utilized, and the proposed computational approaches respect the
group structure.

This paper is organized as follows. In Section \ref{sec:eom}, a Lie
group variational integrator is developed. Optimal control problems
using impulsive controls are studied in Section \ref{sec:oic}, and
optimal control problems with smooth controls are studied in Section
\ref{sec:oc}. Numerical results for a rigid dumbbell spacecraft
 are given in Section \ref{sec:ne}.

\section{Lie group variational integrator on \SE}\label{sec:eom}

The configuration space for the translational and rotational motion
of a rigid body is the special Euclidean group,
$\SE=\Re^3\,\textcircled{s}\,\SO$. We identify the cotangent bundle
$\T^*\SE$ with $\SE\times\se^*$ by left translation, and we identify
$\se^*$ with $\Re^6$ by an isomorphism between $\Re^6$ and $\se$,
and the standard inner product on $\Re^6$. We denote the attitude,
position, angular momentum, and linear momentum of the rigid body by
$(R,x,\Pi,\gamma)\in\T^*\SE$.

The continuous equations of motion are given by
\begin{gather}
\dot{x}=\frac{\gamma}{m},\label{eqn:xdot}\\
\dot{\gamma}=f+u^f,\label{eqn:gamdot}\\
\dot{R}=RS(\Omega),\label{eqn:Rdot}\\
\dot{\Pi}+\Omega\times\Pi = M+u^m,\label{eqn:Pidot}
\end{gather}
where $\Omega\in\Re^3$ is the angular velocity, and
$u^f,u^m\in\Re^3$ are the control force in the inertial frame and
the control moment in the body fixed frame, respectively. The
constant mass of the rigid body is $m\in\Re$, and $J\in\Re^{3\times
3}$ denotes the moment of inertia, i.e. $\Pi=J\Omega$. The map
$S(\cdot):\Re^3\mapsto\so$ is an isomorphism between $\so$ and
$\Re^3$ defined by the condition $S(x)y=x\times y$ for all
$x,y\in\Re^3$.

We assume that the potential is dependent on the position and the
attitude; $U(\cdot):\SE\mapsto\Re$. The corresponding force and the
moment due to the potential are given by
\begin{align}
f & =-\deriv{U}{x},\label{eqn:f}\\
M & =r_{1}\times u_{r_1}+r_{2}\times u_{r_2}+r_{3}\times
u_{r3},\label{eqn:M}
\end{align}
where $r_i,u_{r_i}\in\Re^3$ are the $i$th row vector of $R$ and
$\deriv{U}{R}$, respectively.

Since the dynamics of a rigid body has the structure of a Lagrangian
or Hamiltonian system, they are characterized by symplectic,
momentum and energy preserving properties.  These geometric features
determine the qualitative behavior of the rigid body dynamics, and
they can serve as a basis for theoretical study of rigid body
dynamics.

In contrast, the most common numerical integration methods,
including the widely used Runge-Kutta schemes, neither preserve the
Lie group structure nor these geometric properties. In addition,
standard Runge-Kutta methods fail to capture the energy dissipation
of a controlled system accurately~\cite{jo:marsden}. Additionally,
if we integrate \refeqn{Rdot} by a typical Runge-Kutta scheme, the
quantity $R^T R$ inevitably drifts from the identity matrix as the
simulation time increases. It is often proposed to parameterize
\refeqn{Rdot} by Euler angles or unit quaternions. However, Euler
angles are not global expressions of the attitude since they have
associated singularities. Unit quaternions do not exhibit
singularities, but are constrained to lie on the unit three-sphere
$\S^3$, and general numerical integration methods do not preserve
the unit length constraint. Therefore, quaternions have the same
numerical drift problem. Renormalizing the quaternion vector at each
step tends to break other conservation properties. Furthermore, unit
quaternions, which are diffeomorphic to $\mathrm{SU(2)}$, double
cover $\SO$. So there are inevitable ambiguities in expressing the
attitude.

In~\cite{pro:cca05}, Lie group variational integrators are
introduced by explicitly adapting  Lie group
methods~\cite{ic:iserles} to the discrete variational
principle~\cite{jo:marsden}. They have the desirable property that
they are symplectic and momentum preserving, and they exhibit good
energy behavior for an exponentially long time period. They also
preserve the Euclidian Lie group structure without the use of local
charts, reprojection, or constraints. These geometrically exact
numerical integration methods yield highly efficient and accurate
computational algorithms for rigid body dynamics.   They avoid
singularities and ambiguities.

Using the results presented in~\cite{jo:CMAME05}, a Lie group
variational integrator on $\SE$ for equations
\refeqn{xdot}--\refeqn{Pidot} is given by
\begin{gather}
x_{k+1} = x_k +\frac{h}{m} \gamma_k + \frac{h^2}{2m}\parenth{f_k+u^f_k},\label{eqn:xkp}\\
\gamma_{k+1} = \gamma_k +\frac{h}{2}\parenth{f_k+u^f_k}+\frac{h}{2}\parenth{f_{k+1}+u^f_{k+1}},\label{eqn:gamkp}\\
h S(\Pi_k + \frac{h}{2}\parenth{M_k+u^m_k}) = F_k J_d -J_d F_k^T,\label{eqn:findf}\\
R_{k+1} = R_k F_k,\label{eqn:Rkp}\\
\Pi_{k+1} = F_k^T\Pi_k +\frac{h}{2}F_k^T \parenth{M_k+u^m_k} +
\frac{h}{2}\parenth{M_{k+1}+u^m_{k+1}},\label{eqn:Pikp}
\end{gather}
where the subscript $k$ denotes the $k$th discrete variables for a
fixed integration step size $h\in\Re$. $J_d\in\Re^{3\times 3}$ is a
nonstandard moment of inertia matrix defined by
$J_d=\frac{1}{2}\mathrm{tr}\!\bracket{J}I_{3\times 3}-J$.
$F_k\in\SO$ is the relative attitude between adjacent integration
steps.

For given $(R_k,x_k,\Pi_k,\gamma_k)$ and control inputs,
\refeqn{findf} is solved to find $F_k$. Then $(R_{k+1},x_{k+1})$ are
obtained by \refeqn{Rkp},\refeqn{xkp}. Using \refeqn{f},\refeqn{M},
$(f_{k+1},M_{k+1})$ are computed, and they are used to find
$(\Pi_{k+1},\gamma_{k+1})$ by \refeqn{Pikp},\refeqn{gamkp}. This
yields a map
$(R_k,x_k,\Pi_k,\gamma_k)\mapsto(R_{k+1},x_{k+1},\Pi_{k+1},\gamma_{k+1})$,
and this process is repeated. The only implicit part is
\refeqn{findf}. The actual computation of $F_k$ is done in the Lie
algebra $\so$ of dimension 3, and the rotation matrices are updated
by multiplication. This approach is completely different from
integration of the kinematics equation \refeqn{Rdot}; there is no
excessive computational burden. It can be shown that this integrator
has second order accuracy. The properties of these discrete
equations of motion are discussed in more detail in
\cite{pro:cca05,jo:CMAME05}.

\section{Optimal Impulsive Control of a Rigid Body}\label{sec:oic}
We formulate an optimal impulsive control problem for a rigid body
on $\SE$, and we develop sensitivity derivatives. They are used in
our computational method for solve optimal impulsive control
problems.

\subsection{Problem formulation}
An optimal impulsive control problem is formulated as a maneuver of
a rigid body from a given initial configuration
$(R_0,x_0,\Pi_0,\gamma_0)$ to a desired configuration described by
\begin{align*}
\braces{(R_N,x_N,\Pi_N,\gamma_N)\in\T^*\SE \big|
\mathcal{C}(R_N,x_N,\Pi_N,\gamma_N)=0},
\end{align*}
where $\mathcal{C}(\cdot):\T^*\SE\mapsto\Re^c$ during the given
maneuver time $N$. Two impulsive control inputs are applied at the
initial time and the terminal time. We assume that the control
inputs are purely impulsive, which means that each impulse changes
the momentum of the rigid body instantaneously, but it does not have
any effect on the position and the attitude of the rigid body at
that instant. The motion of the rigid body between the initial time
and the terminal time is uncontrolled. i.e. $u^f_k=u^m_k=0$.  The
performance index is the sum of the magnitudes of the initial
impulse and the terminal impulse. It is equivalent to minimizing the
sums of the initial momentum change and the terminal momentum
change.

We transform this optimal impulsive control problem into a parameter
optimization problem. Let $(\Pi_0^+,\gamma_0^+)$ be the initial
momentum after the initial impulsive control. Then, the terminal
states are determined by the discrete equations of motion, and the
momentum after the terminal impulsive control,
$(\Pi_N^+,\gamma_N^+)$, can be computed by the terminal constraint.
Therefore, the performance index and the constraint are completely
determined by $(\Pi_0^+,\gamma_0^+)$. Thus, the optimal impulsive
control on $\SE$ is formulated as
\begin{gather*}
\text{given}  : (R_0,x_0,\Pi_0,\gamma_0), N\\
\begin{aligned}
\min_{\Pi_0^+,\gamma_0^+} \mathcal{J} & = \norm{\Pi_0^+-\Pi_0}
+\norm{\gamma_0^+-\gamma_0}\vspace*{-0.3cm}\\
& \quad +\norm{\Pi_N^+-\Pi_N}+\norm{\gamma_N^+-\gamma_N},
\end{aligned}\\
\text{such that } \mathcal{C}(R_N,x_N,\Pi_N^+,\gamma_N^+)=0,\\
\text{ subject to discrete equations of motion
\refeqn{xkp}--\refeqn{Pikp}.}
\end{gather*}
If the desired values for all of the terminal states are specified
by the constraints, then there is no freedom for optimization. This
problem degenerates to a two point boundary value problem on $\SE$,
which can be considered as an extension of the Lambert problem for
the restricted two body problem. A similar optimal control problem
for attitude dynamics of a rigid body on $\SO$ is studied
in~\cite{pro:acc06}.

\subsection{Sensitivity derivatives}

\subsubsection*{Variational model}
The variation of $g_k=(R_k,x_k)\in\SE$ can be expressed in terms of
a Lie algebra element $\eta_k\in\se$ and the exponential map as
$g_k^\epsilon = g_k \exp \epsilon\eta_k$. The corresponding
infinitesimal variation is given by
\begin{align*}
\delta g_k & = \frac{d}{d\epsilon}\bigg|_{\epsilon=0}g_k \exp
\epsilon\eta_k = \T_e L_{g_k}\cdot\eta_k.
\end{align*}
Using homogeneous coordinates~\cite{bk:MuLiSa}, the above equation
is written in a matrix equation as
\begin{align}
\ase{\delta R_k}{\delta x_k} & = \aSE{R_k}{x_k}\ase{S(\zeta_k)}{\chi_k},\nonumber\\
& = \ase{R_kS(\zeta_k)}{R_k\chi_k},\label{eqn:infvarg}
\end{align}
where $\zeta_k,\chi_k\in\Re^3$ so that $(S(\zeta_k),\chi_k)\in\se$.
This gives an expression for the infinitesimal variation of a Lie
group element in terms of its Lie algebra. Then, small perturbations
from a given trajectory on $\T^*\SE$ can be written as
\begin{align}
x_k^\epsilon & = x_k + \epsilon\delta x_k,\label{eqn:xke}\\
\gamma_k^\epsilon & = \gamma_k + \epsilon\delta\gamma_k,\\
\Pi_k^\epsilon & = \Pi_k + \epsilon\delta\Pi_k,\\
R_k^{\epsilon} 
& = R_k + \epsilon R_k
S(\zeta_k)+\mathcal{O}(\epsilon^2)\label{eqn:Rke},
\end{align}
where $\delta x_k,\delta \gamma_k,\delta\Pi_k,\zeta_k$ are
considered in $\Re^3$.

We derive expressions for the constrained variation of $F_k$ using
\refeqn{Rkp} and \refeqn{Rke}. Since $F_k=R_k^T R_{k+1}$ by
\refeqn{Rkp}, the infinitesimal variation $\delta F_k$ is given by
\begin{align*}
\delta F_k & = \delta R_k^T R_{k+1} + R_k^T \delta R_{k+1} =
-S(\zeta_k)F_k + F_k S(\zeta_{k+1}).
\end{align*}
We can also express $\delta F_k=F_kS(\xi_k)$ for $\xi_k\in\Re^3$,
using \refeqn{infvarg}. Using the property $S(R^T x)=R^TS(x)R$ for
all $R\in\SO$ and $x\in\Re^3$, we obtain the constrained variation
of $F_k$
\begin{align}
\xi_k = -F_k^T \zeta_k +\zeta_{k+1}.\label{eqn:xik0}
\end{align}

\subsubsection*{Linearized equations of motion} Substituting the
variation model \refeqn{xke}--\refeqn{Rke} and the constrained
variation \refeqn{xik0} into the equations of motion
\refeqn{xkp}--\refeqn{Pikp}, and ignoring higher order terms, the
linearized equation of motion can be written as
\begin{align}
z_{k+1} = A_k z_k,\label{eqn:sl}
\end{align}
where $z_k=[\delta x_{k}; \delta\gamma_{k}; \zeta_{k};
\delta\Pi_{k}]\in\Re^{12}$, and $A_k\in\Re^{12\times 12}$
can be suitably defined. The solution of \refeqn{sl} is obtained as
\begin{align}
z_N 
 = \Phi
z_0,\label{eqn:zN}
\end{align}
where $\Phi\in\Re^{12\times 12}$ represents the sensitivity
derivatives of the terminal state with respect to the initial state
on $\SE$.

\subsection{Computational approach}
We solve the optimal impulsive control problem by the Sequential
Quadratic Programming (SQP) method using analytical expressions for
the gradients of the performance index and the constraints. The
exact computation of the gradients are crucial for efficient
numerical optimization. For the given problem, $\delta
x_0=\zeta_0=0$ since the initial position and the initial attitude
are fixed. Thus, \refeqn{zN} is written as
\begin{align}
\begin{bmatrix}\delta x_{N}\\ \delta\gamma_{N}\\\zeta_{N}\\\delta\Pi_{N}\end{bmatrix}
& = \begin{bmatrix} \Phi^{12} & \Phi^{14}\\
\Phi^{22} & \Phi^{24}\\
\Phi^{32} & \Phi^{34}\\
\Phi^{42} & \Phi^{44}
 \end{bmatrix}
\begin{bmatrix}\delta\gamma_{0}^+\\\delta\Pi_{0}^+\end{bmatrix},\label{eqn:sd}
\end{align}
where $\Phi^{ij}\in\Re^{3\times 3}$, $i,j\in(1,2,3,4)$ are
submatrices of $\Phi$. The above equation represents the
sensitivities of the terminal state with respect to the initial
momentum $(\Pi_0^+,\gamma_0^+)$. Therefore, we can obtain
expressions for gradients of the performance index and the
constraints, and any Newton type numerical approach can be applied.

\section{Optimal Control of a Rigid Body}\label{sec:oc}
We formulate an optimal control problem for a rigid body on $\SE$
assuming that control forces and moments are applied during the
maneuver. Necessary conditions for optimality are developed and
computational approaches are presented to solve the corresponding
two point boundary value problem.

\subsection{Problem formulation}
An optimal impulsive control problem is formulated as a maneuver of
a rigid body from a given initial configuration
$(R_0,x_0,\Pi_0,\gamma_0)$ to a desired configuration
$(R_N^d,x_N^d,\Pi_N^d,\gamma_N^d)$ during the given maneuver time
$N$. Control inputs are parameterized by their value at each time
step. The performance index is the square of the weighted $l_2$ norm
of the control inputs.
\begin{gather*}
\text{given: } (x_0,\gamma_0,R_0,\Pi_0),\,(x_N^d,\gamma_N^d,R_N^d,\Pi_N^d),\,N,\\
\min_{u_{k+1}} \mathcal{J}=\sum_{k=0}^{N-1}
\frac{h}{2}(u^f_{k+1})^TW_fu^f_{k+1}+
\frac{h}{2}(u^m_{k+1})^TW_mu^m_{k+1},\\
\text{such that } (x_N,\gamma_N,R_N,\Pi_N)=(x_N^d,\gamma_N^d,R_N^d,\Pi_N^d),\\
\text{subject to discrete equations of motion
\refeqn{xkp}--\refeqn{Pikp},}
\end{gather*}
where $W_f,W_m\in\Re^{3\times 3}$ are symmetric positive definite
matrices.  Here we use a modified version of the discrete equations
of motion with first order accuracy, because it yields a compact
form for the necessary conditions, which are developed the following
subsection. A similar optimal control problem for attitude dynamics
on $\SO$ is studied in~\cite{JOTA06}.

\subsection{Necessary conditions for optimality} Define an
augmented performance index as
\begin{align*}
\mathcal{J}_a =
\sum_{k=0}^{N-1}&\frac{h}{2}(u^f_{k+1})^TW^fu^f_{k+1}+
\frac{h}{2}(u^m_{k+1})^TW^mu^m_{k+1}\nonumber\\
& +\lambda_k^{1,T}\braces{-x_{k+1}+x_k+\frac{h}{m}\gamma_k}\nonumber\\
& +\lambda_k^{2,T}\braces{-\gamma_{k+1} + \gamma_k+hf_{k+1}+hu^f_{k+1}}\nonumber\\
& +\lambda_k^{3,T}S^{-1}\!\parenth{\mathrm{logm}(F_k-R_{k}^TR_{k+1})}\nonumber\\
& +\lambda_k^{4,T}\braces{-\Pi_{k+1} + F_k^T \Pi_k +
h\parenth{M_{k+1}+u_{k+1}^m}},
\end{align*}
where $\lambda_k^{i}\in \Re^3$ are Lagrange multipliers. The
constraint \refeqn{findf} is considered implicitly using a
constrained variation.
Using the variational model \refeqn{xke}--\refeqn{Rke}, the
constrained variation \refeqn{xik0}, and the fact that the
variations vanish at $k={0,N}$, we obtain the infinitesimal
variation of $\mathcal{J}_a$ as
\begin{align*}
\delta\mathcal{J}_a & = \sum_{k=1}^{N-1} h\delta
u_{k}^{f,T}\braces{W_fu^f_{k}+\lambda_{k-1}^2}\\
&+h\delta u_{k}^{m,T}\braces{W_mu^m_{k}+\lambda_{k-1}^4}
+z_k^T\braces{-\lambda_{k-1}+A_k^T \lambda_k},
\end{align*}
where
$\lambda_k=[\lambda_k^1;\lambda_k^2;\lambda_k^3;\lambda_k^4]\in\Re^{12}$,
and $A_k\in\Re^{12\times 12}$ is presented in \refeqn{sl}.

Since $\delta\mathcal{J}_a=0$ for all variations,  we obtain necessary
conditions for optimality as follows.
\begin{gather}
x_{k+1}=x_k+\frac{h}{m}\gamma_k,\label{eqn:updatex1}\\
\gamma_{k+1} = \gamma_k+hf_{k+1}+h u^f_{k+1},\label{eqn:updategam1}\\
h S(\Pi_k) = F_k J_d - J_dF_k^T,\label{eqn:findf1}\\
R_{k+1} = R_k F_k,\label{eqn:updateR1}\\
\Pi_{k+1} = F_k^T \Pi_k + hM_{k+1}+hu^m_{k+1},\label{eqn:updatePi1}\\
u^f_{k+1} = -W_{f}^{-1}\lambda_{k}^2,\label{eqn:ufkp}\\
u^m_{k+1} = -W_{m}^{-1}\lambda_{k}^4,\label{eqn:umkp}\\
\lambda_{k} = A_{k+1}^T \lambda_{k+1}.\label{eqn:updatelam}
\end{gather}
In the above equations, the only implicit part is \refeqn{findf1}.
For a given initial condition $(R_0,x_0,\Pi_0,\gamma_0)$ and
$\lambda_0$, we can find $F_0$ by solving \refeqn{findf1}. Then,
$R_1,x_1$ is obtained by \refeqn{updateR1},\refeqn{updatex1}, and
the control input $u^f_1,u^m_1$ is obtained by
\refeqn{ufkp},\refeqn{umkp}. $\gamma_1,\Pi_1$ can be obtained by
\refeqn{updategam1},\refeqn{updatePi1}. Now we compute
$(R_1,x_1,\Pi_1,\gamma_1)$. We solve \refeqn{findf1} to find $F_1$.
Finally, $\lambda_1$ can be obtained by \refeqn{updatelam}. This
yields a map $\braces{(R_0,x_0,\Pi_0,\gamma_0),\lambda_0}\mapsto
\braces{(R_1,x_1,\Pi_1,\gamma_1),\lambda_1}$, and this process can
be repeated.

\subsection{Computational Approach}
The necessary conditions for optimality are expressed in terms of a
two point boundary problem on $\T^{*}\SE$ and its dual. This problem
is to find the optimal discrete flow, multiplier, and control inputs
to satisfy the equations of motion
\refeqn{updatex1}--\refeqn{updatePi1}, optimality conditions
\refeqn{ufkp},\refeqn{umkp}, multiplier equations
\refeqn{updatelam}, and boundary conditions simultaneously.

We use a neighboring extremal method~\cite{bk:bryson}. A nominal
solution satisfying all of the necessary conditions except the
boundary conditions is chosen. The unspecified initial multiplier is
updated by successive linearization so as to satisfy the specified
terminal boundary conditions in the limit. This is also referred to
as a shooting method. The main advantage of the neighboring extremal
method is that the number of iteration variables is small. In other
approaches, the initial guess of control input history or multiplier
variables are iterated, so the number of optimization parameters are
proportional to the number of discrete time steps.

The difficulty is that the extremal solutions are sensitive to small
changes in the unspecified initial multiplier values. The
nonlinearities also make it hard to construct an accurate estimate
of sensitivity, and it may result in numerical ill-conditioning.
Therefore, it is important to compute the sensitivities accurately
to apply the neighboring extremal method.

Here the optimality conditions
\refeqn{ufkp} and \refeqn{umkp} are substituted into the equations of motion and the
multiplier equations. The sensitivities of the specified terminal
boundary conditions with respect to the unspecified initial
multiplier conditions is obtained by a linear analysis.

Similar to \refeqn{sl},
the linearized equations of motion can be written as
\begin{align}\label{eqn:zkp}
z_{k+1} & = A_k z_k + \mathcal{A}^{12} \delta\lambda_k,
\end{align}
where
$\mathcal{A}^{12}_k=-h\mathrm{diag}[0,W^{-1}_f,0,W^{-1}_m]\in\Re^{12\times
12}$. We can linearize the multiplier equations \refeqn{updatelam}
to obtain
\begin{align}\label{eqn:dellamkp}
\delta \lambda_k = \mathcal{A}_{k+1}^{21} z_{k+1} + A_{k+1}^T \delta
\lambda_{k+1},
\end{align}
where $\mathcal{A}_{k+1}^{21}\in\Re^{12\times 12}$ can be defined
properly. The solution of the linear equations \refeqn{zkp} and
\refeqn{dellamkp} can be obtained as
\begin{align*}
\begin{bmatrix}z_{N}\\\delta\lambda_{N}\end{bmatrix}
= \begin{bmatrix} \Psi^{11} & \Psi^{12}\\\Psi^{21} &
\Psi^{22}\end{bmatrix}
\begin{bmatrix}z_{0}\\\delta\lambda_{0}\end{bmatrix},
\end{align*}
where $\Psi^{ij}\in\Re^{12\times 12}$.

For the given two point boundary value problem $z_0=0$ since the
initial condition is fixed, and $\lambda_N$ is free.  Thus,
\begin{align}
z_N = \Psi_{12} \delta\lambda_0.
\end{align}
The matrix $\Psi_{12}$ represents the sensitivity of the specified
terminal boundary conditions with respect to the unspecified initial
multipliers. Using this sensitivity, an initial guess of the
unspecified initial conditions is iterated to satisfy the specified
terminal conditions in the limit.

Any type of Newton iteration can be applied. We use a line search
with backtracking algorithm, referred to as Newton-Armijo iteration
in~\cite{bk:kelley}. The procedure is summarized as follows.

\vspace*{0.1cm} {
\renewcommand{\theenumi}{\arabic{enumi}}
\renewcommand{\labelenumi}{\theenumi:}
\newcommand{\tab}{\hspace*{0.6cm}}
\hrule\vspace*{0.08cm}
\begin{enumerate}
\item Guess an initial multiplier $\lambda_0$.
\item Find $x_k,\gamma_k,\Pi_k,R_k,\lambda_k$ using \refeqn{updatex1}--\refeqn{updatelam}.
\item Compute the terminal B.C. error; $\mathrm{Error}=\norm{z_N}$.
\item Set $\mathrm{Error}^t=\mathrm{Error},\;\; i=1$.
\item \textbf{while} $\mathrm{Error} > \epsilon_S$.
\item \tab Find a line search direction; $D=\Psi_{12}^{-1}$.
\item \tab Set $c=1$.
\item \tab\textbf{while} $\mathrm{Error}^t > (1-2\alpha c)\mathrm{Error}$
\item \tab\tab Choose a trial multiplier $\lambda_0^{t}=\lambda_0^{}+c D
z_N$.
\item \tab\tab Find $x_k,\gamma_k,\Pi_k^{},R_k^{},\lambda_k$ using \refeqn{updatex1}--\refeqn{updatelam}.
\item \tab\tab Compute the error; $\mathrm{Error}^t=\norm{z_N^t}$.
\item \tab\tab Set $c=c/10,\;\; i=i+1$.
\item \tab\textbf{end while}
\item \tab Set $\lambda_0=\lambda_0^t$, $\mathrm{Error}=\mathrm{Error}^t$. (accept the trial)
\item \textbf{end while}
\end{enumerate}
\vspace*{0.08cm} \hrule} \vspace*{0.1cm} \noindent Here $i$ is the
number of iterations, and $\epsilon_S,\alpha\in\Re$ are a stopping
criterion and a scaling factor, respectively. The outer loop finds a
search direction by computing the sensitivity derivatives, and the
inner loop performs a line search to find the largest step size
$c\in\Re$ along the search direction. The error in satisfaction of
the terminal boundary condition is determined at each inner
iteration.

\addtolength{\textheight}{-0.12cm}

\section{Numerical Examples}\label{sec:ne}
\subsection{Restricted Full Two Body Problem}
We study a maneuver of a rigid spacecraft under a central gravity
field. We assume that the mass of the spacecraft is negligible
compared to the mass of a central body, and we consider a fixed
frame attached to the central body as an inertial frame. The
resulting model is a Restricted Full Two Body Problem (RF2BP).

The spacecraft is modeled as a dumbbell, which consists of two equal
spheres and a massless rod. The gravitational potential is given by
\begin{align}
U(x,R)=-\frac{GMm}{2} \sum_{q=1}^2 \frac{1}{\norm{x+R\rho^q}},
\end{align}
where $G\in\Re$ is the gravitational constant, $M,m\in\Re$ are the
mass of the central body, and the mass of the dumbbell,
respectively. The vector $\rho^q\in\Re^3$ is the position of the
$q$th sphere from the mass center of the dumbbell expressed in the
body fixed frame ($q\in \braces{1,2}$). The mass, length, and time
dimensions are normalized by the mass of the dumbbell, the radius of
a reference circular orbit, and its orbital period.

\subsection{Optimal Impulsive Control}
We study an impulsive orbital transfer problem with an attitude
change. Initially, the spacecraft is on a reference circular orbit.
We consider two cases. In the first case, the spacecraft moves to a
desired circular orbit and the desired values for all of the
terminal state are specified. There is no freedom for optimization,
and the resulting problem is a two point boundary value problem on
$\SE$. This maneuver can be considered as a generalization of
Hohmann transfer~\cite{bk:danby}. The desired maneuver involves
doubling the orbital radius in addition to a large angle attitude
change.

In the second case, the terminal constraints are relaxed such that
the spacecraft is allowed to transfer to any point on the desired
orbit. The desired terminal orbit is described by its orbital radius
$r_d\in\Re$, and a directional vector $e_n\in\mathbb{S}^2$ normal to
the orbital plane. Two constraints are imposed to locate the
dumbbell in the desired orbital plane with the desired orbital
radius, and one constraint is applied to align the dumbbell to the
normal direction.

\begin{figure}
    \centerline{\subfigure[Spacecraft maneuver]{
    \includegraphics[width=0.57\columnwidth]{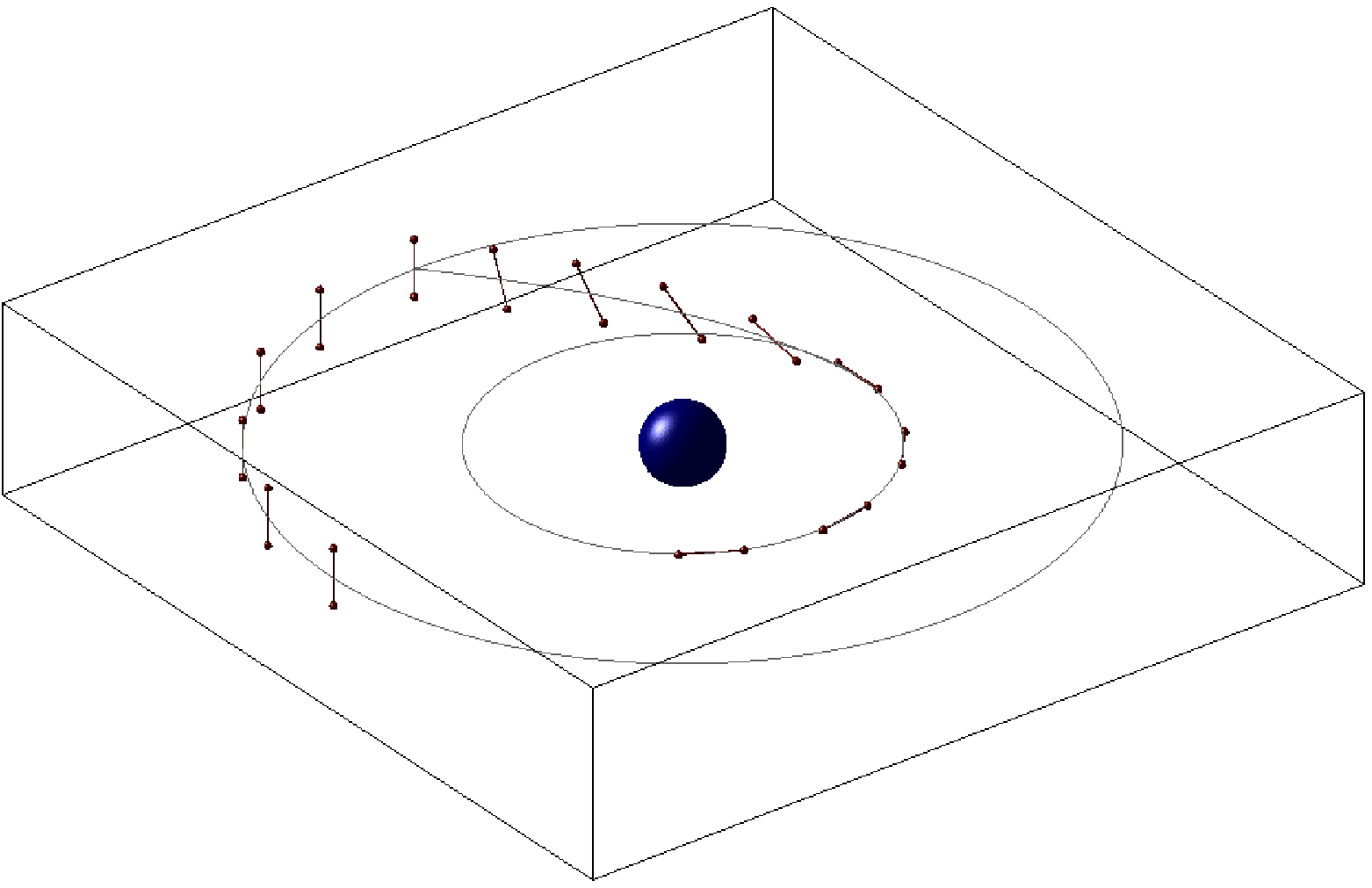}}
    }
    \centerline{\subfigure[Velocity $v$]{
    \includegraphics[width=0.485\columnwidth]{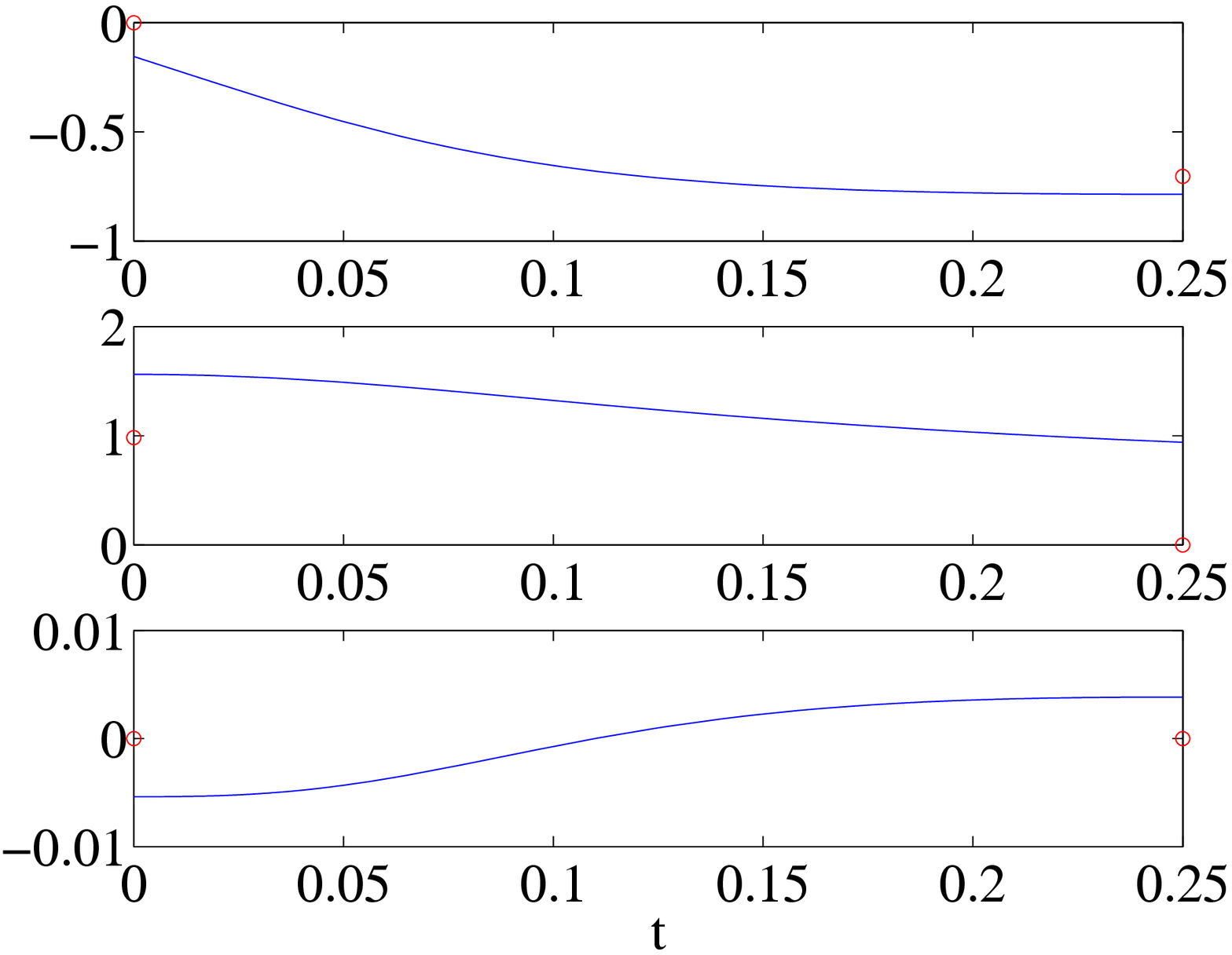}}
    \hfill
    \subfigure[Angular velocity $\Omega$]{
    \includegraphics[width=0.46\columnwidth]{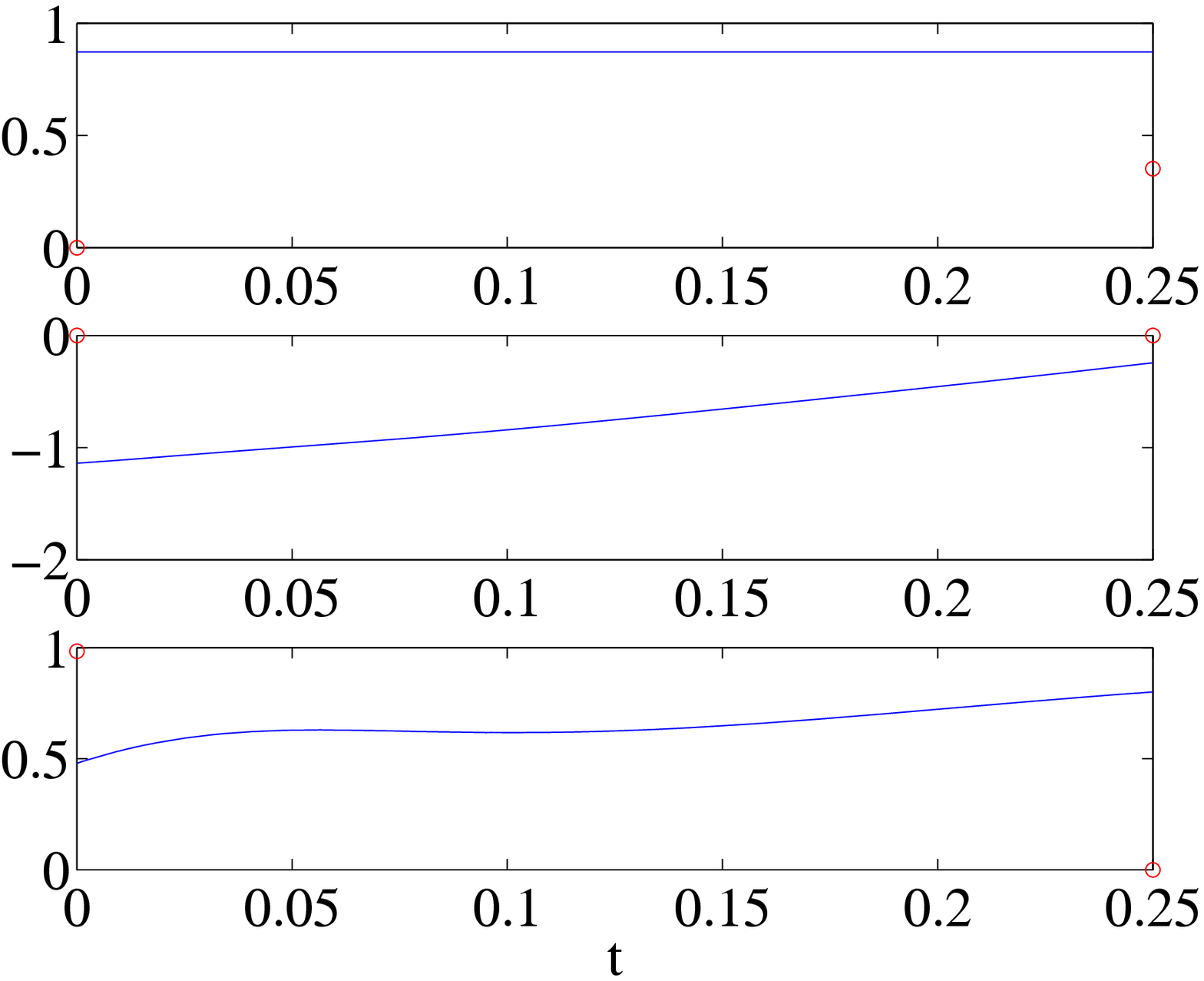}}
    }
    \caption{TPBVP: Orbital radius change}\label{fig:it}
\end{figure}

\begin{figure}
    \centerline{\subfigure[Spacecraft maneuver]{
    \includegraphics[width=0.54\columnwidth]{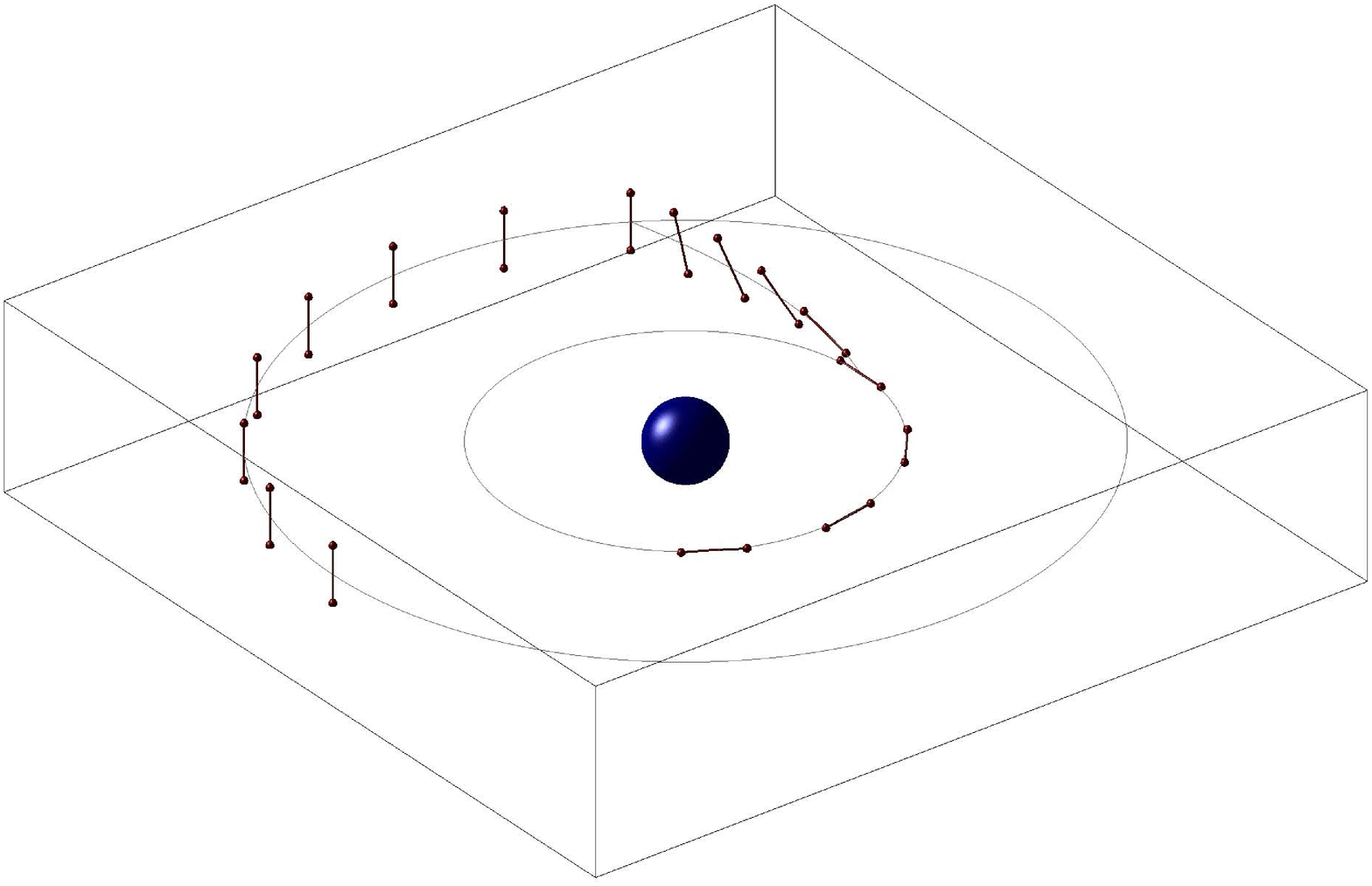}}
    }
    \centerline{\subfigure[Velocity $v$]{
    \includegraphics[width=0.46\columnwidth]{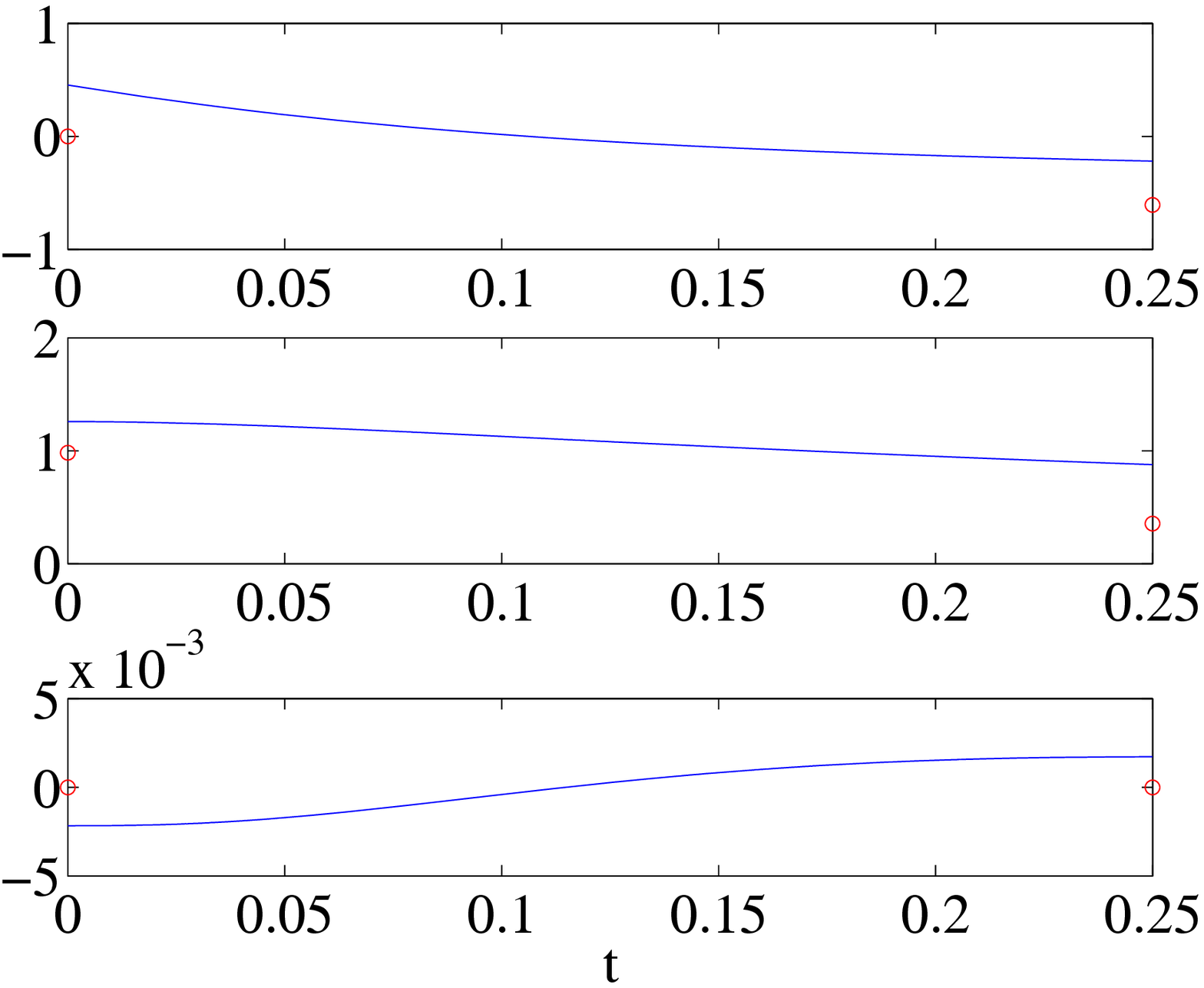}}
    \hfill
    \subfigure[Angular velocity $\Omega$]{
    \includegraphics[width=0.475\columnwidth]{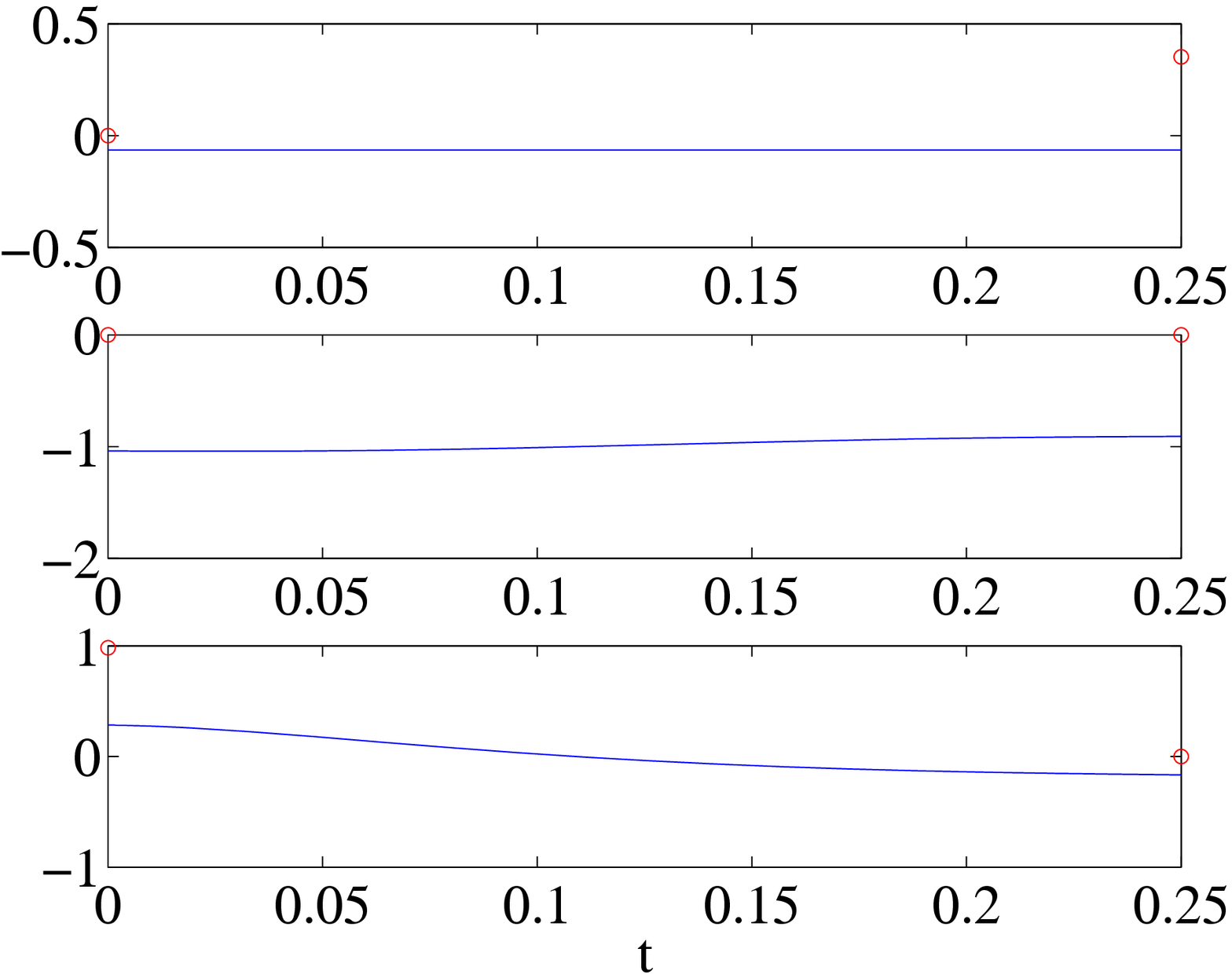}}
    }
    \caption{Optimal impulsive control: Orbital radius change}\label{fig:i}
\end{figure}

The gradients of the performance index and the constraints are
obtained by using \refeqn{sd}. We use Matlab \texttt{fmincon}
function as an implementation of the SQP algorithm. Figures
\ref{fig:it} and \ref{fig:i} show the spacecraft maneuver, and
linear velocity and angular velocity responses, where red circles
denote the velocities before the initial impulse and the velocities
after the terminal impulse. Thus, differences between solid lines
and red circles are proportional to the impulsive controls. (Simple
animations which show these maneuvers of the spacecraft can be found
at \url{http://www.umich.edu/~tylee}.) The error in satisfaction of
the terminal boundary value of the first case is $4.77\times
10^{-15}$. The performance index and the maximum violations of the
constraints for the second case are $1.2305$ and $3.88\times
10^{-15}$, respectively.

\subsection{Optimal Control}
We study an optimal orbital transfer problem to increase the orbital
inclination by $60\,\mathrm{deg}$, and an orbital capture problem to
the reference circular orbit.

\begin{figure}
    \centerline{\subfigure[Spacecraft maneuver]{
    \includegraphics[width=0.39\columnwidth]{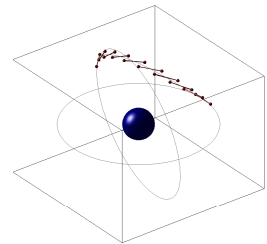}}
    \subfigure[Convergence rate]{
    \includegraphics[width=0.47\columnwidth]{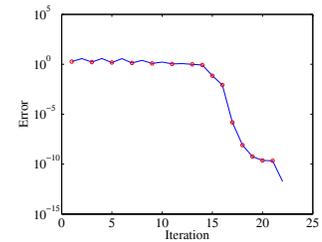}}
    }
    \centerline{\subfigure[Control force $u^f$]{
    \includegraphics[width=0.455\columnwidth]{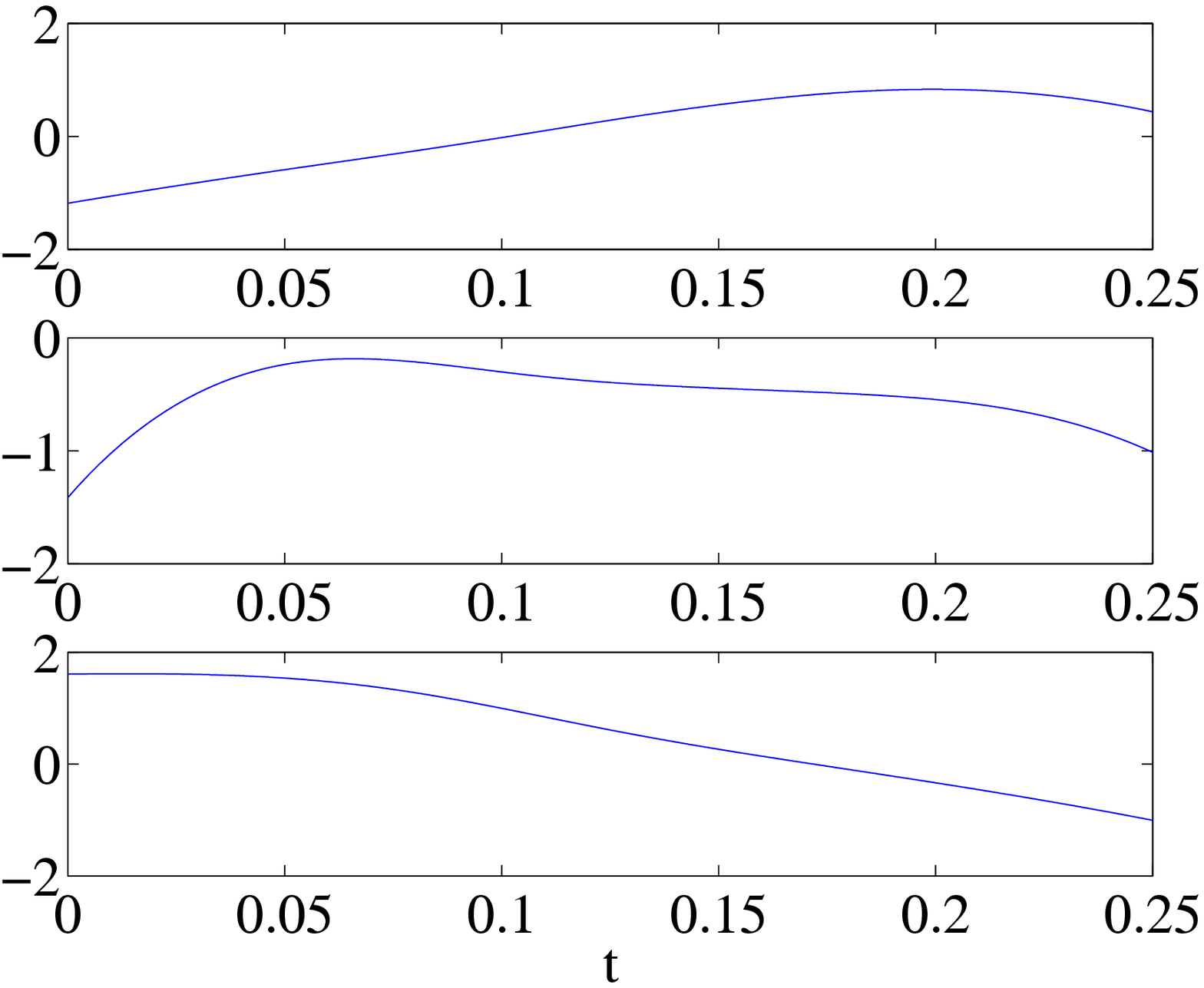}}
    \hfill
    \subfigure[Control moment $u^m$]{
    \includegraphics[width=0.485\columnwidth]{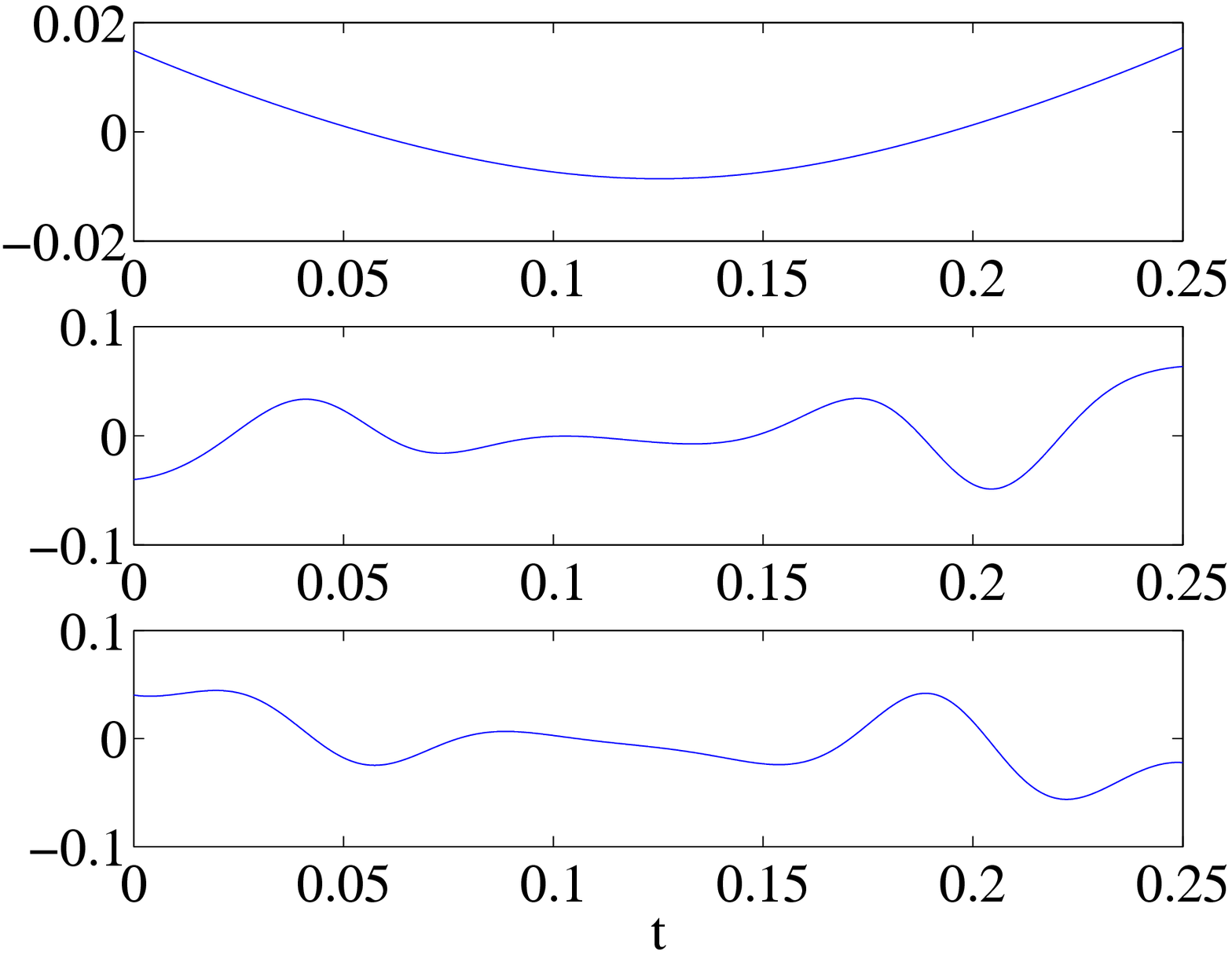}}
    }
    \caption{Optimal control: Orbital inclination change}\label{fig:ii}
\end{figure}

\begin{figure}
    \centerline{\subfigure[Spacecraft maneuver]{
    \includegraphics[width=0.42\columnwidth]{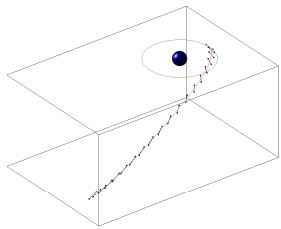}}
    \subfigure[Convergence rate]{
    \includegraphics[width=0.44\columnwidth]{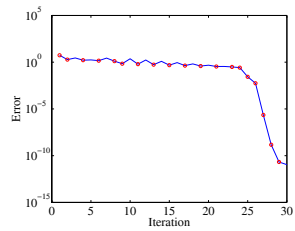}}
    }
    \centerline{\subfigure[Control force $u^f$]{
    \includegraphics[width=0.455\columnwidth]{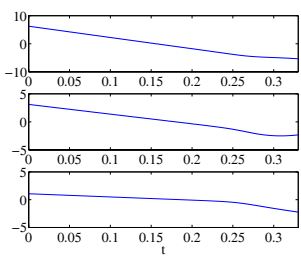}}
    \subfigure[Control moment $u^m$]{
    \includegraphics[width=0.465\columnwidth]{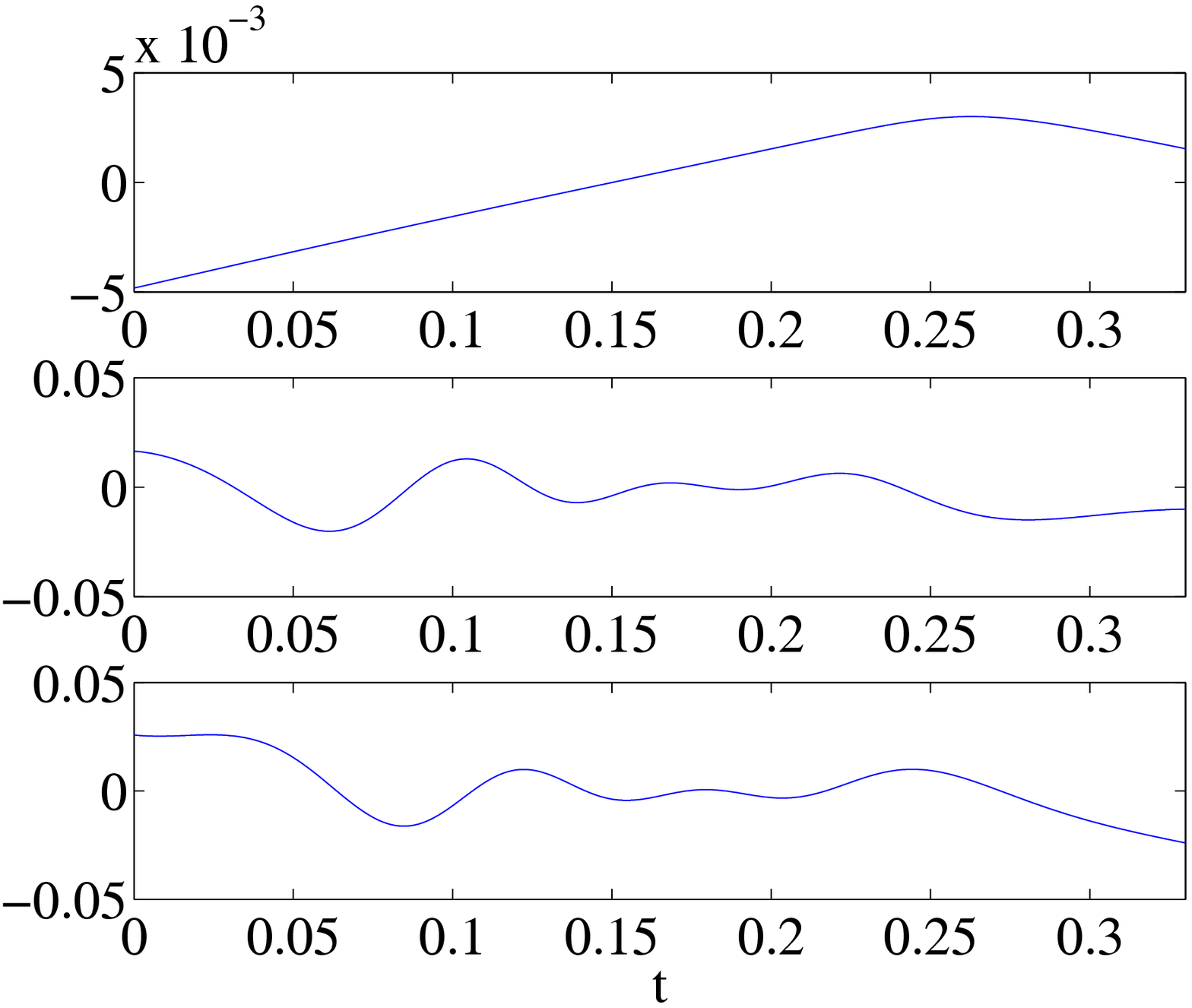}}
    }
    \caption{Optimal control: Orbital capture}\label{fig:iii}
\end{figure}

Figures \ref{fig:ii} and \ref{fig:iii} show the optimized spacecraft
maneuver, control inputs history. For each case, the performance
indices are $13.03$ and $20.90$, and the maximum violations of the
constraint are $3.35\times 10^{-13}$ and $3.26\times 10^{-13}$,
respectively.

Figures \ref{fig:ii}.(b) and \ref{fig:iii}.(b) show the violation of
the terminal boundary condition according to the number of
iterations in a logarithmic scale. Red circles denote outer
iterations in Newton-Armijo iteration to compute the sensitivity
derivatives. For all cases, the initial guesses of the unspecified
initial multiplier are arbitrarily chosen. The error in satisfaction
of the terminal boundary condition converges quickly to machine
precision after the solution is close to the local minimum at around
20th iteration. These convergence results are consistent with the
quadratic convergence rates expected of Newton methods with
accurately computed gradients.

The neighboring extremal method, also referred to as the shooting
method, is numerically efficient in the sense that the number of
optimization parameters is minimized. But, this approach may be
prone to numerical ill-conditioning~\cite{bk:betts}. A small change
in the initial multiplier can cause highly nonlinear behavior of the
terminal attitude and angular momentum. It is difficult to compute
the gradient for Newton iterations accurately, and the numerical
error may not converge.

However, the numerical examples presented in this paper show
excellent numerical convergence properties. This is because the
proposed computational algorithms on $\SE$ are geometrically exact
and numerically accurate.

The dynamics of a rigid body arises from Hamiltonian mechanics,
which have neutral stability, and its adjoint system is also
neutrally stable. The proposed Lie group variational integrator and
the discrete multiplier equations, obtained from variations
expressed in the Lie algebra, preserve the neutral stability
property numerically. Therefore the sensitivity derivatives are
computed accurately.

\section{Conclusions}
Optimal control problems for combined orbital and rotational
maneuvers of a rigid body are formulated and efficient computational
procedures are proposed. The dynamics are discretized by a Lie group
variational integrator, and sensitivity derivatives are developed by
a linear analysis. Discrete necessary conditions for optimality are
constructed, and the corresponding two point boundary value problem
is solved efficiently.

This approach is geometrically exact in the sense that the Lie group
variational integrator preserves the group structure as well as the
geometric invariant properties, and the sensitivity derivatives are
expressed in terms of its Lie algebra. Since the configuration of a
rigid body is defined globally using an element of $\SE$, this
approach completely avoids singularity or ambiguity arising from
other representations such as Euler angles and quaternions.
Numerical examples show the efficiency of the proposed computational
approach.

\addtolength{\textheight}{0.12cm} \clearpage\newpage

\bibliography{opt}

\begin{thebibliography}{10}
\providecommand{\url}[1]{#1}
\csname url@rmstyle\endcsname
\providecommand{\newblock}{\relax}
\providecommand{\bibinfo}[2]{#2}
\providecommand\BIBentrySTDinterwordspacing{\spaceskip=0pt\relax}
\providecommand\BIBentryALTinterwordstretchfactor{4}
\providecommand\BIBentryALTinterwordspacing{\spaceskip=\fontdimen2\font plus
\BIBentryALTinterwordstretchfactor\fontdimen3\font minus
  \fontdimen4\font\relax}
\providecommand\BIBforeignlanguage[2]{{%
\expandafter\ifx\csname l@#1\endcsname\relax
\typeout{** WARNING: IEEEtran.bst: No hyphenation pattern has been}%
\typeout{** loaded for the language `#1'. Using the pattern for}%
\typeout{** the default language instead.}%
\else
\language=\csname l@#1\endcsname
\fi
#2}}

\bibitem{Spin.MCSS98}
K.~Spindler, ``Optimal control on {L}ie groups with applications to attitude
  control,'' \emph{Mathematics of Control, Signals, and Systems}, vol.~11, pp.
  197--219, 1998.

\bibitem{Sas.ICIAM95}
S.~Sastry, ``Optimal control on {L}ie groups,'' in \emph{Proceedings of the
  Third International Congress on Industrial and Applied Mathematics (ICIAM)},
  1995.

\bibitem{pro:cca05}
T.~Lee, M.~Leok, and N.~H. McClamroch, ``A {L}ie group variational integrator
  for the attitude dynamics of a rigid body with applications to the 3{D}
  pendulum,'' in \emph{Proceedings of the IEEE Conference on Control
  Applications}, 2005, pp. 962--967.

\bibitem{jo:marsden}
J.~E. Marsden and M.~West, ``Discrete mechanics and variational integrators,''
  \emph{Acta Numerica}, vol.~10, pp. 357--514, 2001.

\bibitem{ic:iserles}
A.~Iserles, H.~Z. Munthe-Kaas, S.~P. N{\o}rsett, and A.~Zanna, ``Lie-group
  methods,'' \emph{Acta Numerica}, vol.~9, pp. 215--365, 2000.

\bibitem{jo:CMAME05}
\BIBentryALTinterwordspacing
T.~Lee, M.~Leok, and N.~H. McClamroch, ``Lie group variational integrators for
  the full body problem,'' \emph{Computer Methods in Applied Mechanics and
  Engineering}, 2005, submitted. [Online]. Available:
  \url{http://arxiv.org/math.NA/0508365}
\BIBentrySTDinterwordspacing

\bibitem{pro:acc06}
\BIBentryALTinterwordspacing
------, ``Attitude maneuvers of a rigid spacecraft in a circular orbit,'' in
  \emph{Proceedings of the American Control Conference}, 2006, pp. 1742--1747.
  [Online]. Available: \url{http://arxiv.org/math.NA/0509299}
\BIBentrySTDinterwordspacing

\bibitem{bk:MuLiSa}
R.~M. Murray, Z.~Li, and S.~S. Sastry, \emph{A Mathematical Introduction to
  Robotic Manipulation}.\hskip 1em plus 0.5em minus 0.4em\relax CRC Press,
  1993.

\bibitem{JOTA06}
\BIBentryALTinterwordspacing
T.~Lee, M.~Leok, and N.~H. McClamroch, ``Optimal attitude control of a rigid
  body using geometrically exact computations on {S}{O}(3),'' \emph{Journal of
  Optimization Theory and Applications}, 2006, submitted. [Online]. Available:
  \url{http://arxiv.org/math.OC/0601424}
\BIBentrySTDinterwordspacing

\bibitem{bk:bryson}
A.~E. Bryson and Y.-C. Ho, \emph{Applied Optimal Control}.\hskip 1em plus 0.5em
  minus 0.4em\relax Hemisphere Publishing Corporation, 1975.

\bibitem{bk:kelley}
C.~T. Kelley, \emph{Iterative Methods for Linear and Nonlinear
  Equations}.\hskip 1em plus 0.5em minus 0.4em\relax SIAM, 1995.

\bibitem{bk:danby}
J.~M.~A. Danby, \emph{Fundamentals of Celestial Mechanics}.\hskip 1em plus
  0.5em minus 0.4em\relax Willmann Bell Inc., 1988.

\bibitem{bk:betts}
J.~T. Betts, \emph{Practical Methods for Optimal Control Using Nonlinear
  Programming}.\hskip 1em plus 0.5em minus 0.4em\relax SIAM, 2001.

\end{thebibliography}
\bibliographystyle{IEEEtran}

\end{document}